               \newcommand{\be}{\beta}
\newcommand{\ga}{\gamma}               
               
\newcommand{\lb}{\lambda}              \newcommand{\Lb}{\Lambda}
\newcommand{\sig}{\sigma}              
               
\newcommand{\veps}{\varepsilon}        \newcommand{\vphi}{\varphi}

\newcommand{\cal}{\mathcal}
           
           \newcommand{\calf}{{\cal F}}

\newcommand{\card}{{\rm card}}    
\newcommand{\chain}{{\rm chain}}      
      
\newcommand{\Fix}{{\rm Fix}}          
         
\newcommand{\rchain}{{\rm rchain}}

\newcommand{\incl}{\subseteq}        
         
\newcommand{\es}{\emptyset}          
      \newcommand{\limpl}{\Longrightarrow}
  
\newcommand{\oo}{\infty}

             \newcommand{\sk}{\smallskip}

                \def\R+oo{R_+\cup\{\oo\}}

\def\dtends   {\stackrel {\it d}{\longrightarrow}}
\def\0dtends   {\stackrel {\it 0d}{\longrightarrow}}

\newcommand{\barr}{\begin{array}}         \newcommand{\earr}{\end{array}}
\newcommand{\bcor}{\begin{corollary}}     \newcommand{\ecor}{\end{corollary}}
\newcommand{\beq}{\begin{equation}}       \newcommand{\eeq}{\end{equation}}
\newcommand{\bex}{\begin{example}}        \newcommand{\eex}{\end{example}}
\newcommand{\bit}{\begin{itemize}}        \newcommand{\eit}{\end{itemize}}
\newcommand{\blemma}{\begin{lemma}}       \newcommand{\elemma}{\end{lemma}}
\newcommand{\bproof}{\begin{proof}}       \newcommand{\eproof}{\end{proof}}
\newcommand{\bprop}{\begin{proposition}}  \newcommand{\eprop}{\end{proposition}}
\newcommand{\brem}{\begin{remark}}        \newcommand{\erem}{\end{remark}}
\newcommand{\btab}{\begin{tabular}}       \newcommand{\etab}{\end{tabular}}
\newcommand{\btheorem}{\begin{theorem}}   \newcommand{\etheorem}{\end{theorem}}

\documentclass[reqno]{amsart}

\newtheorem{theorem}{\bf Theorem}
\newtheorem{corollary}{\bf Corollary}
\newtheorem{example}{\bf Example}
\newtheorem{lemma}{\bf Lemma}
\newtheorem{proposition}{\bf Proposition}
\newtheorem{remark}{\bf Remark}

\begin{document}

\title
[Functional contractions in local Branciari metric spaces]
{FUNCTIONAL CONTRACTIONS IN \\
LOCAL BRANCIARI METRIC SPACES}

\author{Mihai Turinici}
\address{
"A. Myller" Mathematical Seminar;
"A. I. Cuza" University;
700506 Ia\c{s}i, Romania
}
\email{mturi@uaic.ro}


\subjclass[2010]{
47H10 (Primary), 54H25 (Secondary).
}

\keywords{
Symmetric, polygonal inequality, local Branciari metric,
convergent/Cauchy sequence, Matkowski function, 
contraction, periodic and fixed point.
}

\begin{abstract}
A fixed point result is given for a class of
functional contractions over local Branciari metric spaces. 
It extends some contributions in the area due to
Fora et al [Mat. Vesnik, 61 (2009), 203-208].
\end{abstract}

\maketitle

\section{Introduction}
\setcounter{equation}{0}

Let $X$ be a nonempty set.
By a {\it symmetric} over $X$ we shall mean any map 
$d:X\times X\to R_+$ with 
(cf. Hicks \cite{hicks-1999})
\bit
\item[(a01)]
$d(x,y)=d(y,x)$,\ \ $\forall x,y\in X$ \hfill ($d$ is {\it symmetric});
\eit
the couple $(X,d)$ will be referred to 
as a {\it symmetric space}.
Further, let $T:X\to X$ be a selfmap of $X$,
and put $\Fix(T)=\{z\in X; z=Tz\}$;
any such point will be called {\it fixed} under $T$.
According to Rus \cite[Ch 2, Sect 2.2]{rus-2001},
we say that $x\in X$ is a {\it Picard point} (modulo $(d,T)$) if
{\bf 1a)} $(T^nx; n\ge 0)$ is $d$-convergent, 
{\bf 1b)} each point of $\lim_n(T^nx)$ is in $\Fix(T)$.
If this happens for each $x\in X$
then $T$ is referred to as a {\it Picard operator} (modulo $d$); 
if (in addition) $\Fix(T)$ is singleton ($x,y\in \Fix(T)$ $\limpl$ $x=y$), 
then $T$ is called a {\it strong Picard operator} (modulo $d$). 
[We refer to Section 2 for all unexplained notions].
Sufficient conditions for these properties  to be valid
require some additional conditions upon $d$; 
the usual ones are
\bit
\item[(a02)]
$d(x,y)=0$ iff $x=y$ \ \ \hfill ($d$ is {\it reflexive sufficient})
\item[(a03)]
$d(x,y)\le d(x,z)+d(z,y)$,\ $\forall x,y,z\in X$
\hfill ($d$ is {\it triangular});
\eit
when both these hold, 
$d$ is called a {\it (standard) metric} on $X$.
In this (classical) setting, a basic result 
to the question we deal with is the 1922 one due to 
Banach \cite{banach-1922};
it says that, whenever $(X,d)$ is complete and 
(for some $\lb$ in $[0,1[$)
\bit
\item[(a04)]
$d(Tx,Ty)\le \lb d(x,y)$,\ \ $\forall x,y\in X$,
\eit
then, $T$ is a strong Picard operator.
This result found various applications in 
operator equations theory;
so, it was the subject of many extensions.
A natural way of doing this  
is by considering "functional" contractive conditions
like
\bit
\item[(a05)]
$d(Tx,Ty)\le F(d(x,y),d(x,Tx),d(y,Ty),d(x,Ty),d(y,Tx))$,
\ \ $\forall x,y\in X$;
\eit
where $F:R_+^5\to R_+$ is an appropriate function.
For more details about the possible choices of $F$
we refer to the 1977 paper by
Rhoades \cite{rhoades-1977};
see also 
Turinici \cite{turinici-1976}.
Another way of extension is that of conditions imposed
upon $d$ being modified.
For example, in the class of symmetric spaces,
a relevant paper concerning the contractive question is  
the 2005 one due to
Zhu et al \cite{zhu-cho-kang-2005}.
Here, we shall be interested
in fixed point results 
established over generalized metric spaces,
introduced as in 
Branciari \cite{branciari-2000};
where, the triangular property (a03) is to 
be substituted by the {\it tetrahedral} one:
\bit
\item[(a06)]
$d(x,y)\le d(x,u)+d(u,v)+d(v,y)$,\\
whenever $x,y,u,v\in X$ are distinct to each other.
\eit
Some pioneering results in the area were given by
Das \cite{das-2002}, 
Mihe\c{t} \cite{mihet-2009},
and 
Samet \cite{samet-2009};
see also 
Azam and Arshad \cite{azam-arshad-2008}.
In parallel to such developments, 
some technical problems involving these structures were considered.
For example, Sarma et al \cite{sarma-rao-rao-2009}
observed that Branciari's result may 
not hold, in view of the Hausdorff property for 
$(X,d)$ being not deductible from (a06).
This remark was followed by a series of results 
founded on this property being {\it ab initio} imposed;
see, in this direction,
Chen and Sun  \cite{chen-sun-2012}
or
Lakzian and Samet \cite{lakzian-samet-2011}.
However, in 2011,
Kikina and Kikina \cite{kikina-kikina-2011}
noticed that such a regularity  condition is 
ultimately superfluous for the ambient space;
so, the initial setting will suffice for
these results being retainable.
It is our aim in the present exposition to 
confirm this, within a class of "local"  
Branciari metric spaces.
Further aspects will be delineated elsewhere.

\section{Preliminaries}
\setcounter{equation}{0}

Let $N=\{0,1,...\}$ denote the set of all natural numbers.
For each $n\ge 1$ in $N$, let $N(n,>):=\{0,...,n-1\}$ stand for
the {\it initial interval} (in $N$) induced by $n$.
Any set $P$ with 
$P\sim N$ (in the sense: there exists a bijection from 
$P$ to $N$) will be referred to as 
{\it effectively denumerable}; 
also denoted as: $\card(P)=\aleph_0$.
In addition, given some $n\ge 1$, any set $Q$ with 
$Q\sim N(n,>)$ will be said to be {\it $n$-finite};
and we write this: $\card(Q)=n$;
when $n$ is generic here, we say that $Q$ is {\it finite}.
Finally, the (nonempty) set $Y$ is called
(at most) {\it denumerable} iff it is either 
effectively denumerable or finite.
\sk

{\bf (A)}
Let $(X,d)$ be a symmetric space.
Given $k\ge 1$, any ordered system 
$C=(x_1,...,x_k)$ in $X^k$ will be called a 
{\it $k$-chain} of $X$;
the class of all these will be re-denoted as $\chain(X;k)$.
Given such an object, put $[C]=\{x_1,...,x_k\}$. 
If $\card([C])=k$, then $C$ will be referred to as a 
{\it regular $k$-chain} (in $X$);
denote the class of all these as $\rchain(X;k)$.
In particular, any point $a\in X$ may be identified with
a regular $1$-chain of $X$. 
For any $C\in \chain(X;k)$, where $k\ge 2$, denote 
\bit
\item[]
$\Lb(C)=d(x_1,x_2)+...+d(x_{k-1},x_k)$,\ \ 
whenever $C=(x_1,...,x_k)$
\eit
(the "length" of $C$).
Given $h\ge 1$ and the $h$-chain $D=(y_1,...,y_h)$ in $X$,
let $(C;D)$ stand for the $(k+h)$-chain $E=(z_1,...,z_{k+h})$ in $X$ 
introduced as
\bit
\item[]
$z_i=x_i$, $1\le i\le k$;\ \ $z_{k+j}=y_j$, $1 \le j\le h$;
\eit
it will be referred to as the "product" between $C$ and $D$.
This operation may be extended to 
a finite family of such objects.

Having these precise, let us say that the 
symmetric $d$ is a {\it local Branciari metric} when
it is reflexive sufficient and has the property:
for each effectively denumerable $M\incl X$, 
there exists $k=k(M)\ge 1$ such that
\bit
\item[(b01)]
$d(x,y)\le \Lb(x;C;y)$, for all $x,y\in M$, $x\ne y$, and \\
all $C\in \rchain(M;k)$, with 
$(x;C;y)\in \rchain(M;2+k)$
\eit
(referred to as: the {\it $(2+k)$-polyhedral inequality}).
Note that, the triangular inequality (a03)
and the tetrahedral inequality (a06)
are particular cases of this one,
corresponding to $k=1$ and $k=2$, respectively.
On the other hand, (b01) is not reducible 
to (a03) or (a06); because, 
aside from $k> 2$ being allowed, 
the index in question depends on each 
effectively denumerable subset $M$ of $X$.

Suppose that we introduced such an object.
Define a $d$-convergence structure 
over $X$ as follows.
Given the sequence $(x_n)$ in $X$ and the point $x\in X$,
we say that $(x_n)$, $d$-converges to $x$ 
(written as: $x_n \dtends x$)
provided $d(x_n,x)\to 0$; i.e.,
\bit
\item[(b02)]
$\forall\veps> 0$, $\exists i=i(\veps)$:\ \ 
$n\ge i \limpl d(x_n,x)< \veps$.
\eit
(This concept meets the standard requirements in
Kasahara \cite{kasahara-1976}; 
we do not give details).
The set of all such points $x$ will be denoted $\lim_n (x_n)$;
when it is nonempty, $(x_n)$ is called {\it $d$-convergent}.
Note that, in this last case, 
$\lim_n(x_n)$ may be not a singleton, even if (a06) holds;
cf. Samet \cite{samet-2010}.
Further, call the sequence $(x_n)$, {\it $d$-Cauchy} 
when $d(x_m,x_n)\to 0$ as $m,n\to \oo$, $m< n$; i.e.,
\bit
\item[(b03)]
$\forall\veps> 0$, $\exists j=j(\veps)$:\ \ 
$j\le m< n \limpl d(x_m,x_n)< \veps$.
\eit
Clearly, a necessary condition for this is
\bit
\item[]
$d(x_m,x_{m+i})\to 0$ as $m\to \oo$,\ \ 
for each $i> 0$;
\eit
referred to as: $(x_n)$ is {\it $d$-semi-Cauchy};
but the converse is not in general true. 
Note that, by the adopted setting, 
a $d$-convergent sequence need not be $d$-Cauchy.
even if $d$ is tetrahedral;
see the quoted paper for details.
Despite of this, $(X,d)$ is called {\it complete},
if each $d$-Cauchy sequence is $d$-convergent.
\sk

{\bf (B)}
As already precise, the (nonempty) set of limit points for
a convergent sequence is not in general a singleton.
However, in the usual (metric) fixed point arguments,
the convergence property of this sequence
comes from the $d$-Cauchy property of the same.
So, we may ask whether this supplementary condition
upon $(x_n)$ will suffice for such a property.
Call $(X,d)$, {\it Cauchy-separated}
if, for each $d$-convergent $d$-Cauchy sequence 
$(x_n)$ in $X$, $\lim_n(x_n)$ is a singleton.

\bprop  \label{p1}
Assume that $d$ is a local Branciari metric (see above).
Then, $(X,d)$ is Cauchy-separated.
\eprop

\bproof
Let $(x_n)$ be a $d$-convergent $d$-Cauchy sequence.
Assume by contradiction that $\lim_n(x_n)$ has 
at least two distinct points:
\bit
\item[(b04)]
$\exists u,v\in X$ with $u\ne v$, such that\ 
$x_n\dtends u$, $x_n\dtends v$.
\eit

{\bf i)}
Denote $A=\{n\in N; x_n=u\}$, $B=\{n\in N; x_n=v\}$.
We claim that both $A$ and $B$ are finite.
In fact, if $A$ is effectively denumerable, then
$A=\{n(j); j\ge 0\}$, where $(n(j); j\ge 0)$
is strictly ascending (hence $n(j)\to \oo$ as $j\to \oo$) and 
$x_{n(j)}=u$, $\forall j\ge 0$.
Since, on the other hand, $x_{n(j)}\to v$ as $j\to \oo$, 
we must have $d(u,v)=0$; so that, $u=v$, contradiction.
An identical reasoning is applicable when $B$ is effectively denumerable; 
hence the claim. 
As a consequence, there exists $p\in N$, such that:
$x_n\ne u$, $x_n\ne v$, for all $n\ge p$.
Without loss, one may assume that $p=0$; i.e.,
\beq \label{201}
\mbox{
$\{x_n; n\ge 0\}\cap \{u,v\}=\es$\ \ 
[$x_n\ne u$ and $x_n\ne v$,\ \ for all $n\ge 0$].
}
\eeq

{\bf ii)}
Put $h(0)=0$. We claim that the set $S_0=\{n\in N; x_n=x_{h(0)}\}$
is finite. For, otherwise, it has the representation  
$S_0=\{m(j); j\ge 0\}$, where $(m(j); j\ge 0)$
is strictly ascending (hence $m(j)\to \oo$ as $j\to \oo$) and 
$x_{m(j)}=x_0$, $\forall j\ge 0$.
Combining with (b04) gives $x_0=u$, $x_0=v$; 
hence, $u=v$, contradiction. 
As a consequence of this, there exists $h(1)> h(0)$ 
with $x_{h(1)}\ne x_{h(0)}$.
Further, by a very similar reasoning, 
$S_{0,1}=\{n\in N; x_n\in \{x_{h(0)},x_{h(1)}\}\}$
is finite too; hence, there exists $h(2)> h(1)$ with
$x_{h(2)}\notin \{x_{h(0)},x_{h(1)}\}$; and so on.
By induction, we get a subsequence $(y_n:=x_{h(n)}; n\ge 0)$
of $(x_n)$ with 
\beq \label{202}
\mbox{
$y_i\ne y_j$, for $i\ne j$;\ \ 
$y_n\dtends u$, $y_n\dtends v$ as $n\to \oo$.
}
\eeq
The subset $M=\{y_n; n\ge 0\}\cup\{u,v\}$ 
is effectively denumerable;
let $k=k(M)\ge 1$ stand for the natural number
assured by the local Branciari metric property of $d$.
From the $(2+k)$-polyhedral inequality (b01)
we have, for each $n\ge 0$,
$$
d(u,v)\le d(u,y_{n+1})+...+d(y_{n+k},v).
$$
(The possibility of writing this is assured by 
(\ref{201}) and (\ref{202}) above).
On the other hand, $(y_n)$ is a $d$-Cauchy sequence;
because, so is $(x_n)$; hence $d(y_m,y_{m+1})\to 0$
as $m\to \oo$. Passing to limit in the above 
relation gives $d(u,v)=0$; whence, $u=v$, contradiction.
So, (b04) is not acceptable;  and this concludes the argument.
\eproof

{\bf (B)}
Let $\calf(R_+)$ stand for the class of all functions 
$\vphi:R_+\to R_+$. Denote

\bit
\item[(b05)]
$\calf_r(R_+)=\{\vphi\in \calf(R_+); 
\vphi(0)=0;\ \vphi(t)< t,\ \forall t> 0\}$;
\eit
each  $\vphi\in \calf_r(R_+)$ will be referred to as {\it regressive}.
Note that, for any such function, 
\beq \label{203}
\forall u,v\in R_+:\  v\le\vphi(\max\{u,v\}) \limpl v\le \vphi(u).
\eeq
Call $\vphi\in \calf_r(R_+)$, {\it strongly regressive}, provided
\bit 
\item[(b06)]
$\forall \ga> 0$,\ $\exists \be\in ]0,\ga[$,\ $(\forall t)$:\  
$\ga\le t< \ga+\be \limpl \vphi(t)\le \ga$; \\
or, equivalently:\
$0\le t< \ga+\be \limpl \vphi(t)\le \ga$.
\eit 
Some basic properties of such functions are given below.

\bprop \label{p2}
Let $\vphi\in \calf_r(R_+)$ be strongly regressive. Then,

{\bf i)}
for each sequence $(r_n; n\ge 0)$ in $R_+$ with
$r_{n+1}\le \vphi(r_n)$, $\forall n$, we have $r_n\to 0$
[we then say that $\vphi$ is {\it iteratively asymptotic}]

{\bf ii)}
in addition, for each sequence $(s_n; n\ge 0)$ in $R_+$ 
with $s_{n+1}\le \vphi(\max\{s_n,r_n\})$, $\forall n$
we have $s_n\to 0$.
\eprop

\bproof
i)\
Let $(r_n; n\ge 0)$ be as in the premise of this assertion.
As $\vphi$ is regressive, we have 
$r_{n+1}\le r_n$, $\forall n$.
The sequence $(r_n; n\ge 0)$ is therefore decreasing;
hence $\ga:=\lim_n(r_n)$ exists in $R_+$.
Assume by contradiction that $\ga> 0$; 
and let $\be\in ]0,\ga[$ be the number indicated by the
strong regressiveness of $\vphi$.
As $r_n\ge \ga> 0$ (and $\vphi$=regressive), 
one gets $r_{n+1}< r_n$, $\forall n$; hence, $r_n> \ga$, $\forall n$.
Further, as $r_n\to \ga$, there exists some rank $n(\be)$ 
in such a way that (combining with the above)
$n\ge n(\be)$ $\limpl$ $\ga< r_n< \ga+\be$.
The strong regressiveness of $\vphi$ then gives
(for the same ranks, $n$)
$\ga< r_{n+1}\le \vphi(r_n)\le \ga$; contradiction.
Consequently, $\ga=0$; and we are done.

ii)\
Let $(r_n; n\ge 0)$ and $(s_n; n\ge 0)$ be as in the premise
of these assertions.
Denote for simplicity $t_n:=\max\{s_n,r_n\}$, $n\ge 0$.
As $\vphi$ is regressive, one has (for each $n$)
$r_{n+1}\le r_n\le t_n$, $s_{n+1}\le t_n$; 
hence (for the same ranks) $t_{n+1}\le t_n$.
The sequence $(t_n; n\ge 0)$ is therefore descending; 
wherefrom, $t:=\lim_n(t_n)$ exists in $R_+$
and $t_n\ge t$, $\forall n$.
Assume by contradiction that $t> 0$;
As $r_n\to 0$, there must be some rank $n(t)$ such that
$n\ge n(t)$ $\limpl$ $r_n< t$.
Combining with the above, one gets $t_n\ge t> r_n$, 
for all $n\ge n(t)$; whence $t_n=s_n$, for all $n\ge n(t)$.
But then, the choice of $(s_n; n\ge 0)$ gives
$s_{n+1}\le \vphi(s_n)$, for all $n\ge n(t)$.
This, along with the  first part of the proof, gives
$s_n\to 0$; hence $t_n\to 0$; contradiction.
Consequently, $t_n\to 0$; and, from this, the conclusion follows.
\eproof

Now, let us give two  basic examples of such functions.
 
{\bf B1)} Suppose that $\vphi\in \calf_r(R_+)$ 
is a {\it Boyd-Wong} function \cite{boyd-wong-1969}; i.e.
\bit
\item[(b07)]
$\limsup_{t\to s+}\vphi(t)< s$,\ \ for all $s> 0$.
\eit
Then, $\vphi$ is strongly regressive.
The verification is immediate, by definition;
so, we do not give details.

{\bf B2)}
Suppose that $\vphi\in \calf_r(R_+)$ is
a {\it Matkowski function} \cite{matkowski-1977}; i.e.
\bit
\item[(b08)]
$\vphi$ is increasing;\ \ $\vphi^n(t)\to 0$ as $n\to \oo$,\ 
for all $t> 0$.
\eit
[Here, for each $n\ge 0$, $\vphi^n$ stands for the $n$-th
iterate of $\vphi$].
Then, $\vphi$ is strongly regressive.
The verification of this assertion is 
to be found to Jachymski \cite{jachymski-1994};
however, for completeness reasons, we shall 
provide it, with some modifications.
Assume by contradiction that 
$\vphi$ is not strongly regressive;
that is (for some $\ga> 0$)
\bit
\item[]
$\forall \be\in ]0,\ga[$, $\exists t\in [\ga,\ga+\be[$: $\vphi(t)> \ga$
(hence, $\ga < t< \ga+\be$).
\eit
As $\vphi$=increasing, this yields $\vphi(t)> \ga$, $\forall t> \ga$.
By induction, we get $\vphi^n(t)> \ga$, for all $n$, and all $t> \ga$. 
Fixing some $t> \ga$, we have
(passing to limit as $n\to \oo$) $0\ge \ga$, contradiction; 
hence the claim.

\section{Main result}
\setcounter{equation}{0}

Let $X$ be a nonempty set; 
and $d(.,.)$ be a local Branciari metric over it, with
\bit
\item[(c01)]
$(X,d)$ is complete (each $d$-Cauchy sequence is $d$-convergent).
\eit
Note that, by Proposition \ref{p1},
for each $d$-Cauchy sequence $(x_n)$ in $X$, 
$\lim_n(x_n)$ is a (nonempty) singleton, $\{z\}$;
as usually, we write $\lim_n(x_n)=\{z\}$  as $\lim_n(x_n)=z$.

Let $T:X\to X$ be a selfmap of $X$.
We say that $x\in X$ is a {\it Picard point} (modulo $(d,T)$) if
{\bf 3a)} $(T^nx; n\ge 0)$ is $d$-Cauchy (hence $d$-convergent), 
{\bf ii)} $\lim_n(T^nx)$ is in $\Fix(T)$.
If this happens for each $x\in X$
then $T$ is referred to as a {\it Picard operator} (modulo $d$); 
if (in addition) $\Fix(T)$ is a singleton, 
then $T$ is called a {\it strong Picard operator} (modulo $d$). 

Now, concrete circumstances guaranteeing such properties
involve functional contractive (modulo $d$) conditions upon $T$.
Precisely, denote for $x,y\in X$:
\bit
\item[(c02)]
$M(x,y)=\max\{d(x,y),d(x,Tx),d(y,Ty)\}$.
\eit
It is easy to see that
\beq \label{301}
M(x,Tx)=\max\{d(x,Tx),d(Tx,T^2x)\},\ \ \forall x,y\in X.
\eeq
Given $\vphi\in \calf_r(R_+)$,
we say that $T$ is {\it $(d,M;\vphi)$-contractive} if
\bit
\item[(c03)]
$d(Tx,Ty)\le \vphi(M(x,y))$, $\forall x,y\in X$.
\eit
The main result of this note is

\btheorem \label{t1}
Suppose that $T$ is $(d,M;\vphi)$-contractive, 
where $\vphi\in \calf_r(R_+)$ is strongly regressive.
Then, $T$ is a strong Picard operator (modulo $d$).
\etheorem

\bproof
First, we check the singleton property.
Let $z_1,z_2\in \Fix(T)$ be arbitrary fixed.
By this very choice,
$$
M(z_1,z_2)=\max\{d(z_1,z_2),0,0\}=d(z_1,z_2).
$$
Combining with the contractive condition yields
$$
d(z_1,z_2)=d(Tz_1,Tz_2)\le \vphi(d(z_1,z_2));
$$
wherefrom $d(z_1,z_2)=0$; hence $z_1=z_2$; 
so that, $\Fix(T)$ is (at most) a singleton.
It remains now to establish the Picard property.
Fix some $x_0\in X$;
and put $x_n=T^nx_0$, $n\ge 0$.
There are several steps to be passed.

{\bf I)}
If $x_n=x_{n+1}$ for some $n\ge 0$, we are done.
So, it remains to discuss the remaining situation;
i.e. (by the reflexive sufficiency of $d$)
\bit
\item[(c04)]
$\rho_n:=d(x_n,x_{n+1})> 0$,\ \ for all $n\ge 0$.
\eit
By the contractive property and (\ref{301}),
$\rho_{n+1}\le \vphi(\max\{\rho_n,\rho_{n+1}\})$, for all $n\ge 0$;
so that (taking (\ref{203}) into account)
\beq \label{302}
\rho_{n+1}\le \vphi(\rho_n),\ \ \forall n\ge 0. 
\eeq
Combining with (c04) one gets that $(\rho_n; n\ge 0)$ is strictly descending; 
moreover, by Proposition \ref{p2}, $\rho_n\to 0$ as $n\to \oo$.

{\bf II)}
Fix $i\ge 1$, and put $(\sig^i_n:=d(x_n,x_{n+i})$; $n\ge 0)$.
Again by the contractive condition, we get the evaluation
$$
\sig^i_{n+1}=d(Tx_n,Tx_{n+i})\le \vphi(M(x_n,x_{n+i}))=
\vphi(\max\{\sig^i_n,\rho_n,\rho_{n+i}\}),\ \forall n\ge 0;
$$
wherefrom, by (\ref{302}) 
\beq \label{303}
\sig^i_{n+1}\le \vphi(\max\{\sig^i_n,\rho_n\}),\ \ \forall n\ge 0.
\eeq
This yields (again via Proposition \ref{p2})
$\sig^i_n\to 0$, for each $i\ge 1$; that is,
\beq \label{304}
\mbox{
$d(x_n,x_{n+i})\to 0$ as $n\to \oo$, for each $i\ge 1$;
}
\eeq
or, in other words: $(x_n)$ is $d$-semi-Cauchy.

{\bf III)}
Suppose that
\bit
\item[(c05)]
there exists $i,j\in N$ such that $i< j$, $x_i=x_j$.
\eit
Denoting $p=j-i$, we thus have $p> 0$ and $x_i=x_{i+p}$;
so that (by the very definition of our iterative sequence)
$$
\mbox{
$x_i=x_{i+np}$, $x_{i+1}=x_{i+np+1}$,\ \  for all $n\ge 0$.  
}
$$
By the introduced notations this yields (via (c04) and (\ref{304}))
$$
\mbox{
$0< \rho_i=\rho_{i+np}\to 0$ as $n\to \oo$;
}
$$
contradiction. 
Hence, (c05) cannot hold; wherefrom, we must have 
\beq  \label{305}
\mbox{
for all $i,j\in N$:\ \ $i\ne j$ implies $x_i\ne x_j$.
}
\eeq

{\bf IV)}
As a consequence of this fact, the map $n\mapsto x_n$ is injective;
so that, $Y:=\{x_n; n\ge 0\}$ is effectively denumerable.
Let $k=k(Y)\ge 1$ be the natural number attached to it, by the local 
Branciari property of $d$. Also, let $\ga> 0$ be arbitrary fixed;
and $\be\in ]0,\ga[$ be given by 
the strong regressivity of $\vphi$.
By the $d$-semi-Cauchy property (\ref{304}), 
there exists $j(\be)\in N$ such that
\beq \label{306}
d(x_n,x_{n+i})< \be/4k\ (< \be/2<\ga+\be/2),\ \ 
\forall n\ge j(\be),\ \forall i\in \{1,...,k+1\}.
\eeq
We now claim that
\beq \label{307}
(\forall p\ge 1):\ 
d(x_n,x_{n+p})< \ga+\be/2,\ \  \forall n\ge j(\be);
\eeq
and, from this, the $d$-Cauchy property for $(x_n; n\ge 0)$ follows.
The case of $p\in \{1,...,k+1\}$ is clear, via (\ref{306}).
Assume that (\ref{307}) holds, for some $p\ge k+1$;
we show that it holds as well for $p+1$.
So, let $n\ge j(\be)$ be arbitrary fixed.
By the inductive hypothesis and (\ref{306}),
$$ \barr{l}
d(x_{n+k},x_{n+p})< \ga+\be/2< \ga+\be \\
d(x_{n+k},x_{n+k+1})< \be/4k< \be< \ga+\be \\
d(x_{n+p},x_{n+p+1})< \be/4k< \be< \ga+\be;
\earr
$$
whence, by definition,
$$
M(x_{n+k},x_{n+p})< \ga+\be.
$$
This, by the contractive condition and (b06), gives
$$
d(x_{n+k+1},x_{n+p+1})\le 
\vphi(M(x_{n+k},x_{n+p}))\le \ga.
$$
Combining with the $(2+k)$-polyhedral inequality
(applied to $C=(x_{n+2},...,x_{n+k+1})$)
$$  \barr{l}
d(x_n,x_{n+p+1})\le 
d(x_n,x_{n+2})+...+d(x_{n+k},x_{n+k+1})+d(x_{n+k+1},x_{n+p+1}) \\
< \be(k+1)/4k+\ga \le  \be/2 +\ga;
\earr
$$
and the assertion follows.
As $(X,d)$ is complete, we have
\beq \label{308}
\mbox{
$x_n\dtends x$ as $n\to \oo$, for some $z\in X$;
}
\eeq
moreover, by Proposition \ref{p1}, $z$ is uniquely
determined by this relation.
We claim that this is our desired point.
Assume by contradiction that $z\ne Tz$; 
or, equivalently, $\rho:=d(z,Tz)> 0$.

{\bf V)}
Denote 
$A=\{n\in N; x_n=z\}$, $B=\{n\in N; x_n=Tz\}$.
If $A$ is effectively denumerable, we have
$A=\{m(j); j\ge 0\}$, where $(m(j); j\ge 0)$ is
strictly ascending (hence $m(j)\to \oo$). 
As $x_{m(j)}=z$, $\forall j\ge 0$, 
we have $x_{m(j)+1}=Tz$, $\forall j\ge 0$.
Combining with $x_{m(j)+1}\dtends z$ as $j\to \oo$,
we must have $d(z,Tz)=0$; hence $z=Tz$, contradiction.
On the other hand, if $B$ is effectively denumerable, 
we have 
$B=\{n(j); j\ge 0\}$, where $(n(j); j\ge 0)$ is
strictly ascending (hence $n(j)\to \oo$). 
As $x_{n(j)}=Tz$, $\forall j\ge 0$,
one gets (via $x_{n(j)}\dtends z$ as $j\to \oo$) 
$d(z,Tz)=0$; whence $z=Tz$, again a contradiction.
It remains to discuss the case of 
both $A$ and $B$ being finite; i.e.,
\bit
\item[(c06)]
there exists $h\ge 0$ such that: 
$\{x_n; n\ge h\}\cap\{z,Tz\}=\es$.
\eit
The subset $Y:=\{x_n; n\ge h\}\cup \{z,Tz\}$ 
is therefore effectively denumerable.
Let $k=k(Y)\ge 1$ be the natural number
attached to it, by  the local Branciari property of $d$.
We have, for each $n\ge k$ 
(by the $(2+k)$-polyhedral inequality applied to 
$C:=(x_{n+2},...,x_{n+k+1})$)
\beq \label{309}  \
\barr{l}
\rho\le 
d(z,x_{n+2})+...+d(x_{n+k},x_{n+k+1})+d(x_{n+k+1},Tz) \\
\earr
\eeq
By (\ref{304}) and (\ref{308}), there exists $j(\rho)\ge h$ 
in such a way that 
$$
n\ge j(\rho) \limpl d(x_n,z), d(x_n,x_{n+1})< \rho/2.
$$
As a consequence, we must have
$$
M(x_{n+k},z)=\rho,\ \ \forall n\ge j(\rho).
$$
so that, by  the contractive condition,
$$
d(x_{n+k+1},Tz)\le \vphi(\rho), \forall n\ge j(\rho).
$$
Replacing in (\ref{309}), we get an evaluation like
$$ \barr{l}
\rho\le 
d(z,x_{n+2})+...+d(x_{n+k},x_{n+k+1})+\vphi(\rho), \ \ 
\forall n\ge j(\rho).
\earr
$$
Passing to limit as $n$ tends to infinity gives
$\rho\le \vphi(\rho)$; wherefrom (as $\vphi$ is regressive)
$\rho=0$; contradiction.
Hence,  $z=Tz$; and the proof is complete.
\eproof

In particular, when the regressive function $\vphi$
is a Boyd-Wong one,
our main result covers the one due to
Das and Dey \cite{das-dey-2009};
note that, by the developments in 
Jachymski \cite{jachymski-2011},
it includes as well the related statements in 
Di Bari and Vetro \cite{di-bari-vetro-2011}.
On the other hand, when $d(.,.)$ is a standard metric,
Theorem \ref{t1} reduces to the statement in
Leader \cite{leader-1979}.
Further aspects may be found in 
Kikina et al \cite{kikina-kikina-gjino-2012};
see also
Khojasteh et al \cite{khojasteh-razani-moradi-2011}.

\section{Further aspects}
\setcounter{equation}{0}

A direct inspection of the proof above shows that
conclusion of Theorem \ref{t1} is retainable even
if one works with orbital completeness of the ambient
space. Some conventions are in order.
Let $X$ be a nonempty set; and $d(.,.)$ be 
a reflexive sufficient symmetric over it;
supposed to be a local Branciari metric.
Further, take a selfmap $T$ of $X$.
Call the sequence $(y_n; n\ge 0)$ in $X$, {\it $T$-orbital}
when $y_n=T^nx$, $n\ge 0$, for some $x\in X$.
In this case, let us say that $(X,d)$ is 
{\it $T$-orbital complete} when each $T$-orbital  
$d$-Cauchy sequence is $d$-convergent.

The following extension of Theorem \ref{t1} 
is available. Let the general conditions above be 
fulfilled; as well as
\bit
\item[(d01)]
$(X,d)$ is $T$-orbital complete.
\eit

\btheorem \label{t2}
Suppose that $T$ is $(d,M;\vphi)$-contractive, 
where $\vphi\in \calf_r(R_+)$ is strongly regressive.
Then, $T$ is a strong Picard operator (modulo $d$).
\etheorem

The proof mimics the one of Theorem \ref{t1};
so, we omit it.

Call the regressive function $\vphi\in \calf(R_+)$, 
{\it admissible} provided
\bit
\item[(d02)]
$\vphi$ is increasing, usc and 
$\sum_n\vphi^n(t)< \oo$,\ $\forall t> 0$.
\eit
Clearly, $\vphi$ is a Matkowski function; 
hence, in particular, a strongly regressive one.
This, in the particular case of $d$
fulfilling the tetrahedral inequality,
tells us that the main result in 
Fora et al \cite{fora-bellour-al-bsoul-2009}
is a particular case of Theorem \ref{t2} above.
In addition, we note that the usc condition 
posed by the authors may be removed.
Note that, the introduced framework 
may be also used to get an extension of 
the contributions due to
Akram and Siddiqui \cite{akram-siddiqui-2003};
see also
Moradi and Alimohammadi
\cite{moradi-alimohammadi-2011}.
These will be discussed elsewhere.


\end{document}